\magnification=1200
\baselineskip=0.2in
%
%
%
%
\def\a{{\alpha}}

\def\s{{\sigma}}
\def\D{{\Delta}}
\def\Z{{\bf Z}}
\def\e{{\epsilon}}
\def\d{{\delta}}
\def\a{{\alpha}}

\def\A{{\cal A}}

\def\t{{\theta}}

\def\mapright#1{\smash{\mathop{\longrightarrow}\limits^{#1}}}

\def\mapemb#1{\smash{\mathop{\hookrightarrow}\limits^{#1}}}
\font\mysmall=cmr8 at 8pt
\centerline{\bf ON ALGEBRAIC SUPERGROUPS AND}
\centerline {\bf QUANTUM DEFORMATIONS}
\bigskip
\centerline{R. FIORESI\footnote*{Investigation supported by
the University of Bologna, funds for selected research topics.}}
\bigskip
\centerline{Dipartimento di Matematica, Universita' di Bologna}
\centerline{Piazza Porta San Donato 5, 40126 Bologna, Italy}
\centerline{{\mysmall e-mail: fioresi@dm.unibo.it}}
\bigskip
{\bf Abstract. {\mysmall We give the definitions of affine algebraic
supervariety and affine algebraic supergroup 
through the functor of points and we relate them to the
other definitions present in the literature. We study in detail the
algebraic supergroups $GL(m|n)$ and $SL(m|n)$ 
and give explicitly the Hopf algebra
structure of the algebra representing 
the functors of points. 
At the end we give also the quantization of $GL(m|n)$
together with its coaction
on suitable quantum spaces according to Manin's philosophy.
}}
\bigskip
{\bf 1. Introduction}
\medskip
The mathematical foundations of supergeometry
were laid in the 60s by Berezin in [Be] and later by Leites [Le]
Kostant [Ko] and Manin [Ma1]
among many others, its origins being  mainly tied up with 
physical problems.

A new attention to the subject came later 
with the study of quantum
fields and superstring. In the 1999 
``Notes on Supersymmetry'' Deligne and Morgan [DM]
give a categorical point of view on supersymmetry 
notions developed originally by physicists and known from a more 
``operational'' point of view.

In the current
definitions of supermanifold, the points of a supermanifold
are points of an usual manifold and the adjective
super refers to an additional structure on the structural
sheaf of functions on the manifold.
This sheaf is assumed to be a sheaf of commutative superalgebras, where a 
superalgebra is a $\Z_2$-graded algebra. 
When dealing with algebraic supergroups however, there are in some
sense {\it true} points. In fact, for 
each supercommutative algebra $A$ the $A$ points of a supergroup can be viewed
as a certain subset of automorphisms of a superspace $A^{m|n}$
where $A$ is a given commutative superalgebra.
For this reason the point of view of functor of points 
is in this case most useful.
It allows to associate to an affine supergroup a commutative
Hopf superalgebra the same way it does in the commutative case.
It is hence possible to define a quantum deformation
of an algebraic supergroup in complete analogy to
the non super case: it will be a deformation of the
commutative Hopf superalgebra associated to it.

The quantum supergroup $GL(m|n)$ was first constructed by
Manin [Ma3] together with its coactions on suitable quantum
superspaces. The construction of its Hopf superalgebra
structure is in some sense implicit. In a subsequent work, Lyubashenko
and Sudbery [LS] provided quantum deformations of supergroups 
of $GL(m|n)$ type using the universal $R$-matrix formalism.
The same formalism is also used in [P] where more explicit
formulas are given.
In both works however the ring and its coalgebra structure do not
appear explicitly, since the calculation was too involved
using the $R$-matrix approach.

In the present paper we present a definition
of affine algebraic supervariety and
supergroup that is basically equivalent to the one of Manin [Ma1]
and then we use this point of view to give
a quantum deformation of the supergroup $GL(m|n)$.

The main result of the paper consists in giving explicitly
the coalgebra structures for the Hopf superalgebras associated to
the supergroup $GL(m|n)$ and its quantization $k_q[GL(m|n)]$ obtained
according to the Manin phylosophy, that is together with 
coactions on suitable quantum spaces.
The explicit forms of the comultiplication for
$GL(m|n)$ and its quantization $k_q[GL(m|n)]$ are non trivial,
they rely heavily on the presence of the nilpotents
of the superalgebras and do not appear in any other work.
These results were announced in the proceeding [Fi] where
they appeared without proof.


The organization of this paper is as follows.

In \S 2 we introduce the notion of affine supervariety and affine supergroup
using the functor of points. These two definitions turn out to be
basically equivalent to the definitions that one finds in the
literature ([De], [Ma1] among many others).

In \S 3 we write explicitly the Hopf algebra
structure of the Hopf superalgebra associated to 
the supergroups $GL(m|n)$ and $SL(m|n)$, where with $GL(m|n)$ we
intend the supergroup whose $A$ points are the group of automorphisms 
of the superspace $A^{m|n}$ and with $SL(m|n)$
the subsupergroup of $GL(m|n)$ of automorphisms with Berezinian equal to 1,
with $A$ commutative superalgebra.

In \S 4 we construct the non commutative Hopf superalgebra 
$k_q[GL(m|n)]$ deformation in the quantum group sense
of the Hopf superalgebra associated
to $GL(m|n)$. 
We also see that $k_q[GL(m|n)]$ admits coactions on suitable
quantum superspaces.
\medskip
The author wishes to thank Prof. V. S. Varadarajan, Prof. I. Dimitrov,
Prof. A. Vistoli, Prof. A. Sudbery,
Prof. T. Lenagan and Prof. M. A. Lledo for helpful comments.
\bigskip
{\bf 2. Preliminaries on algebraic supervarieties and supergroups}
\medskip

Let $k$ be an algebraically closed field.
All algebras and superalgebras
have to be intended over $k$ unless otherwise specified.
Given a superalgebra $A$ we will denote with $A_0$ the even part,
with $A_1$ the odd part and with $I^A_{odd}$ the ideal generated
by the odd part. 

A superalgebra is said to be {\it commutative} (or {\it supercommutative})
if 
$$
xy=(-1)^{p(x)p(y)}yx, \qquad \hbox{for all homogeneous } x,y
$$
where $p$ denotes the parity of an homogeneous element
($p(x)=0$ if $x \in A_0$, $p(x)=1$ if $x \in A_1$). 

In this section all superalgebras are assumed to be commutative.

Let's denote with $\A$ the category of {\it affine superalgebras}
that is commutative superalgebras such that, modulo the ideal
generated by their odd part, they are affine algebras (an affine
algebra is a finitely generated reduced commutative algebra).

\vskip1ex
{\bf Definition (2.1)}. Define {\it affine algebraic supervariety} 
over $k$ a representable functor $V$ from the category
${\cal A}$ of affine superalgebras 
to the category ${\cal S}$ of sets. 
Let's call $k[V]$ the commutative
$k$-superalgebra representing the functor $V$,
$$
V(A)=Hom_{k-superalg}(k[V],A), \qquad A \in \A.
$$
We will call $V(A)$ the {\it $A$-points} of  the variety $V$.

A morphism of affine supervarieties is identified with a morphism
between the representing objects, that is a morphism of affine
superalgebras.

We also define the functor $V_{red}$ associated to $V$ from the 
category $\A_c$ of affine $k$-algebras to the
category of sets:
$$
V_{red}(A_c)=Hom_{k-alg}(k[V]/I^{k[V]}_{odd},A_c), \qquad A_c \in \A_c.
$$
$V_{red}$ is an affine algebraic variety and it is called the {\it reduced
variety} associated to $V$.

If the algebra $k[V]$ representing the functor $V$
has the additional structure of a commutative
Hopf superalgebra, we say that $V$ is an {\it affine algebraic
supergroup.} (For the definition and main properties of Hopf algebras
see [Mo]).
\medskip
{\bf Remarks (2.2)}. 

1. Let $G$ be an affine algebraic supergroup in the sense of (2.1).
As in the classical setting, the condition
$k[G]$ being a commutative Hopf superalgebra makes the functor
group valued, that is the product
of two morphisms is still a morphism. 

In fact let $A$ be a commutative superalgebra and let
$x, y \in Hom_{k-superalg}(k[G],A)$ be two points of $G(A)$.
The product of $x$ and $y$ is defined as:
$$
x \cdot y=_{def}m_A \cdot x \otimes y \cdot \D
$$
where $m_A$ is the multiplication in $A$ and $\D$ the comultiplication
in $k[G]$. One can directly check that
$x \cdot y \in Hom_{k-superalg}(k[G],A)$, that is:
$$
(x \cdot y)(ab)=(x \cdot y)(a)(x \cdot y)(b)
$$
This is an important
difference with the quantum case that will be treated in 
\S 4. The non commutativity of the Hopf 
algebra in the quantum setting 
does not allow to multiply morphisms(=points).
In fact in the quantum (super)group setting the product of two morphisms
is not in general a morphism.
For more details see [Ma4] pg 13.
\medskip
2. Let $V$ be an affine algebraic supervariety as defined in (2.1).
Let $k_0 \subset k$ be a subfield of $k$. We say that $V$
is a $k_0$-variety if there exists a 
$k_0$-superalgebra $k_0[V]$ such that
$k[V] \cong k_0[V] \otimes_{k_0}k$ and
$$
V(A)=Hom_{k_0-superalg}(k_0[V],A)=Hom_{k-superalg}(k[V],A),
\qquad A \in \A.
$$
We obtain a functor that we still denote by $V$ from the
category $\A_{k_0}$ of affine \break
$k_0$-superalgebras to the category of sets:
$$
V(A_{k_0})=Hom_{k_0-superalg}(k_0[V],A_{k_0}),
\qquad A \in \A_{k_0}.
$$

This allows to consider rationality questions on a supervariety.
We will not pursue this further in the present work.
\vskip1ex
{\bf Examples (2.3).}

1. The $k$-points of an affine supervariety $V$ correspond to the affine
variety defined over $k$ whose functor of points is $V_{red}$.

2. Let $A$ be a commutative superalgebra. 
Let $M(m|n)(A)$ be the linear endomorphisms of the superspace 
$A^{m|n}$ (see [De] pg. 53): 
$$
\pmatrix{a_{11} & \dots  & a_{1,m} & \a_{1,m+1} & \dots & \a_{1,m+n} \cr
         \vdots &        & \vdots & \vdots & & \vdots \cr
         a_{m,1} & \dots  & a_{m,m} & \a_{m,m+1} & \dots & \a_{m,m+n} \cr
       \a_{m+1,1} & \dots  & \a_{m+1,m} & a_{m+1,m+1} & \dots & a_{m+1,m+n} \cr
        \vdots &        & \vdots & \vdots & & \vdots \cr
       \a_{m+n,1} & \dots  & \a_{m+n,m} & a_{m+n,m+1} & \dots & a_{m+n,m+n} \cr
}
$$
$a_{ij} \in A_0$, $\a_{kl} \in A_1$,   
$1 \leq i,j \leq m$ or  $m+1 \leq i,j \leq m+n$, 
$1 \leq k \leq m$, $m+1 \leq l \leq m+n$ or 
$m+1 \leq k \leq m+n$, $1 \leq l \leq m$.
\medskip
This is an affine supervariety represented by the commutative
superalgebra: \break $k[M(m|n)]=$ $k[x_{ij},\xi_{kl}]$
where $x_{ij}$'s 
and $\xi_{kl}$'s are respectively even and odd variables 
with
$1 \leq i,j \leq m$ or  $m+1 \leq i,j \leq m+n$,
$1 \leq k \leq m$, $m+1 \leq l \leq m+n$ or
$m+1 \leq k \leq m+n$, $1 \leq l \leq m$.

Observe that $M_{red}=M(m) \times M(n)$
where $M(l)$ is the functor corresponding to the
affine variety of $l \times l$ matrices.

\medskip
We now would like to give an equivalent point of view and define
again the category of affine supervarieties and affine supergroups.
(See [Ma1], [De], [Be]).
\medskip
{\bf Definition (2.4)}. Let $V_{red}$ be an affine algebraic
variety defined over $k$ and $O_{V_{red}}$ the structural
sheaf of $V_{red}$. Define {\it affine algebraic supervariety} 
$V$ the couple $(V_{red}, O_V)$
where $O_V$ is a sheaf of affine superalgebras such that its stalk is local
and $O_V/I_V$ is isomorphic to $O_{V_{red}}$, where $I_V$ is the 
sheaf of ideals generated by the nilpotent elements.
\medskip
We want to show that this definition is equivalent to the one
given previously.
Clearly if we have an affine supervariety according to the definition
(2.4), we have a superalgebra associated to it, namely the global sections
of the sheaf $O_V$. This means that we 
immediately have the functor of
points associated to it, hence a supervariety according to the definition
(2.1). 

Conversely assume we have a functor of points $V$ 
and a commutative superalgebra
$k[V]$ to which it is associated (see definition (2.1)).
We need to show that it 
gives rise to a sheaf of superalgebras on the affine variety $V_{red}$.

Let's look at the maps:
$$
k[V]_0 \mapemb{\alpha} k[V] \mapright{\beta} k[V]/I^{k[V]}_{odd}
$$
where $I^{k[V]}_{odd}$ is the ideal generated by the nilpotent elements in $k[V]$.
Observe that we have a surjective map $\gamma=\beta \cdot \alpha$.
whose kernel consists of the nilpotent elements of $k[V]_0$.
This induces a map
$$
{\hbox{Spec}} k[V]/I^{k[V]}_{odd} \longrightarrow {\hbox{Spec}} k[V]_0$$
that is an isomorphism since the
kernel of $\gamma$ consists of nilpotents.
Let's now view $k[V]$ as an $k[V]_0$-module. This allows us to build
a sheaf on ${\hbox{Spec}} k[V]_0={\hbox{Spec}} k[V]/I^{k[V]}_{odd}$ 
of $k[V]_0$-modules, where the stalk coincides
with the localization of $k[V]$ into the maximal ideals of $k[V]_0$.
So we have obtained from a commutative superalgebra $k[V]$ a sheaf
of superalgebras on ${\hbox{Spec}} k[V]/I^{k[V]}_{odd}$
which corresponds to the affine variety $V_{red}$, 
whose global sections coincide
with $k[V]$. 
\medskip
The next section will be devoted to construct in detail examples of affine
supergroups.
\bigskip
{\bf 3. The affine supergroups $GL(m|n)$ and $SL(m|n)$}
\medskip
In this section we intend to give explicitly the supergroup
structure for the supergroup functors $GL(m|n)$ and $SL(m|n)$.
For any $A \in \A$ let's define $GL(m|n)(A)$ as the group of
automorphisms of the superspace $A^{m|n}$ 
(see [De] pg 59). 
Define also $SL(m|n)(A)$ as the subset of $GL(m|n)$
of automorphisms with berezinian equal to 1. The berezinian
of the matrix:
$$
\pmatrix{G_{11} & \Gamma_{12} \cr
        \Gamma_{21} & G_{22} } \in SL(m|n)(A), 
$$
($G_{11}$, $G_{22}$ are $m \times m$, $n \times n$ invertible matrices
of even elements, $\Gamma_{12}$, $\Gamma_{21}$ 
are $m \times n$, $n \times m$ matrices
of odd elements) 
is defined as:
$$
Ber=det(G_{22})^{-1} det(G_{11}-\Gamma_{12}G_{22}^{-1}\Gamma_{21}).
$$
(See [Be] ch. 4, [De] pg 59 for more details).

\medskip
From section 2 we know that we need
to give Hopf commutative superalgebras \break $k[GL(m|n)]$ and
$k[SL(m|n)]$
such that for each commutative superalgebra $A$,
$$
\matrix{
GL(m|n)(A)=Hom_{k-superalg}(k[GL(m|n)],A) \cr\cr
SL(m|n)(A)=Hom_{k-superalg}(k[SL(m|n)],A).
}
$$

We start by defining the supercommutative algebras $k[GL(m|n)]$
and $k[SL(m|n)]$, then we will give explicitly its coalgebra and Hopf
algebra structure. 
\vskip1ex

Let $x_{ij}$ for $1 \leq i,j \leq m$ or 
                         $m+1 \leq i,j \leq m+n$ be even variables 
and    $\xi_{kl}$ for $1 \leq k \leq m, m+1 \leq l \leq m+n$ or
                         $m+1 \leq k \leq m+n, 1 \leq l \leq m$
be odd variables.

Denote by $X_{11}$, $X_{22}$, $\Xi_{12}$, $\Xi_{21}$
the following matrices of indeterminates:
$$
\matrix{
X_{11}=(x_{ij})_{1 \leq i,j \leq m}, \qquad
X_{22}=(x_{ij})_{m+1 \leq i,j \leq m+n}, \cr\cr
\Xi_{12}=(\xi_{kl})_{1 \leq k \leq m, m+1 \leq l \leq m+n},\qquad
\Xi_{21}=(\xi_{kl})_{m+1 \leq k \leq m+n, 1 \leq l \leq m}.
}
$$
%
\medskip
{\bf Definition (3.1).}
$$
\matrix{
k[GL(m|n)]=_{def}
k[x_{ij},\xi_{kl},{d^{1...m}_{1...m}}^{-1},
{d^{m+1...m+n}_{m+1...m+n}}^{-1}]
\cr\cr
k[SL(m|n)]=_{def}
{k[x_{ij},\xi_{kl},{d^{1...m}_{1...m}}^{-1},
{d^{m+1...m+n}_{m+1...m+n}}^{-1}]
\over
\pmatrix{
det(S_{22}(X_{22}))
det(X_{11}-\Xi_{12}S_{22}(X_{22})\Xi_{21})-1
}
}
}
$$
where ${d^{1...m}_{1...m}}^{-1}$, ${d^{m+1...m+n}_{m+1...m+n}}^{-1}$
are even variables
such that 
$$
{d^{1...m}_{1...m}}^{-1}d^{1...m}_{1...m}=1, \qquad
{d^{m+1...m+n}_{m+1...m+n}}^{-1}{d^{m+1...m+n}_{m+1...m+n}}=1.
$$
with $d^{1...m}_{1...m}=det(X_{11})$,
$d^{m+1...m+n}_{m+1...m+n}=det(X_{22})$.

To simplify the notation we will also write 
$({d^{1...m}_{1...m}}^{-1})^t$ as ${d^{1...m}_{1...m}}^{-t}$ for
$t$ a positive integer.

$$
\matrix{
S_{11}(x_{ij})=_{def}(-1)^{i-j}A_{ji}^{11}{d^{1...m}_{1...m}}^{-1}, 
\qquad 1 \leq i,j \leq m
\cr \cr
S_{22}(x_{ij})=_{def}(-1)^{i-j}A_{ji}^{22}{d^{m+1...m+n}_{m+1...m+n}}^{-1}, 
\qquad m+1 \leq i,j \leq m+n
}
$$
where $A_{ji}^{11}$ and $A_{ji}^{22}$ denote the 
determinants of the minors obtained by
suppressing the $j^{th}$ row and $i^{th}$ column in $X_{11}$
and $X_{22}$ respectively.

Regarding $S_{11}(X_{11})$ as a matrix of indeterminates
that is $(S_{11}(X_{11}))_{ij}=S_{11}(x_{ij})$,
$1 \leq i,j \leq m$,
we have that:
$$
S_{11}(X_{11})X_{11}=X_{11}S_{11}(X_{11})=I_m
$$
Similarly:
$$
S_{22}(X_{22})X_{22}=X_{22}S_{22}(X_{22})=I_n
$$
where $I_m$ and $I_n$ denote the identity matrix of order $m$ and
$n$ respectively.

The expression: 
$$
Ber=det(S_{22}(X_{22}))det(X_{11}-\Xi_{12}S_{22}(X_{22})\Xi_{21})
$$ 
is called
the {\it Berezinian function}. (See [Be] ch. 3 for more details).
\medskip
{\bf Proposition (3.2)}. 
{\it In the ring $k[GL(m|n)]$ the Berezinian function is invertible}.

{\bf Proof}. Let's write $Ber$ as:
$$
Ber=det(S_{22}(X_{22})det(I_m-\Xi_{12}S_{22}(X_{22})\Xi_{21}
S_{11}(X_{11}))det(X_{11}).
$$
It is enough to prove that $det(I_m-\Xi_{12}S_{22}(X_{22})\Xi_{21}
S_{11}(X_{11}))$ is invertible. To simplify the notation let's call
$A=\Xi_{12}S_{22}(X_{22})$, $B=\Xi_{21}S_{11}(X_{11})$.
We now prove that $I_m-AB$ is invertible (as matrix of indeterminates).
It's inverse is given by:
$$
I_m+(AB)+(AB)^2+ \dots +(AB)^{mn+1}
$$
In fact if one multiplies this matrix with $I_m-AB$,
one obtains a telescopic sum with last term $(AB)^{m+n}$.
Each element of $A$ and $B$ is of degree 1 in the odd indeterminates
hence $(AB)^{mn+1}$ is of odd degree $2mn+2$ hence zero.

\medskip
We now proceed to give the coalgebra structure.
We need to give the comultiplication $\D$ and counit 
$\e$ for all the generators
and verify that they are well defined.

In order make the formulas more readable we introduce the notation: 
$$
a_{ij}=_{def}\cases{x_{ij} & if $1 \leq i,j \leq m$ or 
                         $m+1 \leq i,j \leq m+n$ \cr \cr
              \xi_{ij} & if $1 \leq i \leq m, m+1 \leq j \leq m+n$ or
                         $m+1 \leq i \leq m+n, 1 \leq j \leq m$}
$$
\vskip1ex

{\bf Observation (3.3)}. Let's compute 
$\D(d_{1 \dots m}^{1 \dots m})$ and
$\D(d_{m+1 \dots m+n}^{m+1 \dots m+n})$.
$$
\matrix{
\D(d_{1 \dots m}^{1 \dots m})=
\sum_{\s \in S_m}(-1)^{l(\s)}\D(a_{1,\s(1)}) \dots \D(a_{m,\s(m)})=
\cr\cr
=\sum_{\s \in S_m,1 \leq k_1, \dots , k_m \leq m+n}
(-1)^{l(\s)}a_{1,k_1} \dots a_{m,k_m} \otimes 
a_{k_1,\s(1)} \dots a_{k_m,\s(m)}= \sum_{i=1}^m r_1
}
$$  
with
$$
\matrix{
r_i=\sum_{\s \in S_m,(k_1, \dots , k_m)\in \rho_i}
(-1)^{l(\s)}a_{1,k_1} \dots a_{m,k_m} \otimes
a_{k_1,\s(1)} \dots a_{k_m,\s(m)} \cr\cr
\rho_i=\{(k_1, \dots , k_m) | 
1 \leq k_1, \dots , k_m \leq m+n, \quad
|\{k_1, \dots , k_m\} \cap \{1 \dots m\}|=m-i\}
}
$$
Similarly one can write:
$$
\D(d_{m+1 \dots m+n}^{m+1 \dots m+n})=\sum_{i=1}^n s_i
$$
with
$$
\matrix{
s_i=\sum_{\t \in S_n,(l_1, \dots , l_n)\in \s_i}
(-1)^{l(\t)}a_{1,l_1} \dots a_{n,l_n} \otimes
a_{l_1,\t(1)} \dots a_{l_n,\t(n)} \cr\cr
\s_i=\{(l_1, \dots , l_n) | 1 \leq l_1, \dots , l_n \leq m+n,
\quad
|\{l_1, \dots , l_n\} \cap \{m+1 \dots m+n\}|=n-i\}.
}
$$
\medskip

{\bf Proposition (3.4).}
{\it $k[GL(m|n)]$ and $k[SL(m|n)]$ are bialgebras with comultiplication:
$$
\matrix{
\D(a_{ij})=\sum a_{ik} \otimes a_{kl} \cr\cr
\D({d^{1...m}_{1...m}}^{-1})=
\sum_{i=1}^{2mn+2} (-1)^{i-1}{d^{1...m}_{1...m}}^{-i}
\otimes {d^{1...m}_{1...m}}^{-i} 
(\D(d^{1...m}_{1...m})-{d^{1...m}_{1...m}} \otimes d^{1...m}_{1...m})^{i-1}
\cr \cr \cr
\D({d^{m+1...m+n}_{m+1...m+n}}^{-1})=
\sum_{i=1}^{2mn+2}(-1)^{i-1} {d^{m+1...m+n}_{m+1...m+n}}^{-i}
\otimes {d^{m+1...m+n}_{m+1...m+n}}^{-i}
\cr \cr
(\D(d^{m+1...m+n}_{m+1...m+n})-{d^{m+1...m+n}_{m+1...m+n}} 
\otimes d^{m+1...m+n}_{m+1...m+n})^{i-1}
}
$$
and counit:
$$
\matrix{
\e(a_{ij})=\d_{ij}, \qquad
\e({d^{1...m}_{1...m}}^{-1})=1, \qquad \e({d^{m+1...m+n}_{m+1...m+n}}^{-1})=1.
}
$$}

{\bf Proof}.
One can directly check that these maps are well defined with respect to
the commutation relations among even and odd elements.

For $k[SL(m|n)]$:
$Ber$ is a grouplike element, that is $\D(Ber)=Ber \otimes Ber$.
A proof of this in the quantum case is available in [LS], [P], it
clearly applies also to the non quantum case, that is when $q=1$.

Hence we have immediately:
$$
\D(Ber-1)=(Ber-1) \otimes Ber + 1 \otimes (Ber-1)
$$
We now need to check (for both algebras the computation is
the same) that:
$$
\matrix{
(\D({d^{1...m}_{1...m}}^{-1}))(\D(d^{1...m}_{1...m}))= 1 \otimes 1\cr \cr
(\D({d^{m+1...m+n}_{m+1...m+n}}^{-1}))
(\D(d^{m+1...m+n}_{m+1...m+n}))= 1 \otimes 1
}
$$

By observation (3.3) we have that:
$$
\matrix{
\D(d^{1...m}_{1...m})=r_0+ \dots +r_m,
\quad
\D(d^{m+1...m+n}_{m+1...m+n})=s_0+ \dots +s_n,
\cr\cr
\D({d^{1...m}_{1...m}}^{-1})=
\sum_{i=1}^{2mn+2}(-1)^{i-1}r_0^{-i}(r_1 + \dots + r_m)^{i-1}
\cr\cr
\D({d^{m+1...m+n}_{m+1...m+n}}^{-1})=
\sum_{i=1}^{2mn+2}(-1)^{i-1}s_0^{-i}(s_1 + \dots + s_n)^{i-1}
}
$$
Notice that $r_0=d^{1...m}_{1...m} \otimes d^{1...m}_{1...m}$
and $s_0=d^{m+1...m+n}_{m+1...m+n} \otimes d^{m+1...m+n}_{m+1...m+n}$.

So:
$$
\matrix{
(\D(d^{1...m}_{1...m}))(\D({d^{1...m}_{1...m}}^{-1}))= 
[r_0+ (r_1 +\dots +r_m)][r_0^{-1}-r_0^{-2}(r_1+\dots +r_m)+ \dots ]
\cr\cr
1 \otimes 1 + (r_1 +\dots r_m)r_0-r_0(r_1 +\dots r_m)
-r_0^{-2}(r_1+\dots +r_m)^2+\dots \qquad \qquad (*)
}
$$
We obtain a telescopic sum. 
The generic term is given by:
$$
\matrix{
g=[r_0+(r_1 +\dots +r_m)]
[\dots +(-1)^{i-1}r_0^{-i}(r_1 + \dots + r_m)^{i-1}+
\cr\cr
+(-1)^{i}r_0^{-i-1}(r_1 + \dots + r_m)^{i}+ \dots]=
\cr\cr
=\dots +(-1)^{i-1}r_0^{-i+1}(r_1 + \dots + r_m)^{i-1}+
(-1)^{i-1}r_0^{-i}(r_1 + \dots + r_m)^i+\cr\cr
+(-1)^{i}r_0^{-i}(r_1 + \dots + r_m)^{i}
+(-1)^{i}r_0^{-i+1}(r_1 + \dots + r_m)^{i+1}+ \dots
}
$$
Notice that the second and third term cancel out.

The last term in the sum $(*)$ is given by:
$$
\matrix{
(-1)^{2mn+1}r_0^{-2mn-2}(r_1 + \dots + r_m)^{2mn+1}=
\cr\cr 
(-1)^{2mn+1}r_0^{-2mn-2}
\sum_{1 \leq i_1, \dots,  i_{2mn+1} \leq m}r_{i_1} \dots r_{i_{2mn+1}}.
}
$$
But 
$$
r_{i_1} \dots r_{i_{2mn+1}}=0
$$
since it contains the product of $2mn+1$ odd indeterminates
(each $r_i$ for $i > 1$ contains at least one odd indeterminate).
Hence we have the result.

The check:
$$
\D(d^{m+1...m+n}_{m+1...m+n})\D({d^{m+1...m+n}_{m+1...m+n}}^{-1})=
1 \otimes 1
$$ 
is done in the same way.

Finally one directly checks that $\e$ is a counit.
\medskip
{\bf Remarks (3.5)}. 

1. For $n=m=1$ we have:
$$
\matrix{
\D({(d^1_1)}^{-1})=\D({x_{11}}^{-1})=
{x_{11}}^{-1} \otimes {x_{11}}^{-1}-
{x_{11}}^{-2} \otimes {x_{11}}^{-2}(\xi_{12} \otimes \xi_{21})
\cr\cr
\D({(d^2_2)}^{-1})=\D({x_{22}}^{-1})=
{x_{22}}^{-1} \otimes {x_{22}}^{-1}-
{x_{22}}^{-2} \otimes {x_{22}}^{-2}(\xi_{21} \otimes \xi_{12})
}
$$
that gives precisely the formulas in example (2.3)(3).

2. Notice that if one replaces the $\xi_{kl}$ 
with commuting coordinates, $k[GL(m|n)]$ 
and $k[SL(m|n)]$ are not Hopf algebras. 
This comes from the fact that
the product of two $m+n$ by $m+n$ 
matrices whose diagonal $m \times m$ and $n \times n$ blocks
are invertible is not a matrix of the same type.
The coalgebra and Hopf algebra structures use in
an essential way the supercommutativity (i.e. the presence
of nilpotents).     
\medskip
Let's now define the antipode $S$.

Let $B=(b_{ij})_{m+1 \leq i \leq m+n, 1 \leq j \leq m}$, 
$C=(c_{kl})_{1 \leq k \leq m, m+1 \leq l \leq m+n}$ 
be the following matrices of even elements:
$$
\matrix{
B=X_{11}-\Xi_{12}S_{22}(X_{22})\Xi_{21} \cr \cr
C=X_{22}-\Xi_{21}S_{11}(X_{11})\Xi_{12}
}
$$
Define the matrices:
$$
\matrix{
S_1(B)_{ij}=S_1(b_{ij})=_{def}
(-1)^{i-j}A^B_{ji}det(X_{11}-\Xi_{12}S_{22}(X_{22})\Xi_{21})^{-1}
\cr \cr
S_2(C)_{kl}=S_2(c_{kl})=_{def}
(-1)^{l-k}A^C_{lk}det(X_{22}-\Xi_{21}S_{11}(X_{11})\Xi_{12})^{-1}
}
$$
where $A^B_{ji}$ and $A^C_{lk}$ are the determinants of the minors obtained
by suppressing the $j^{th}$ row and $i^{th}$ column in $B$ and
the $l^{th}$ row and $k^{th}$ column in $C$ respectively.
\medskip
{\bf Remark (3.6)}.
The determinants that appear in the definition of $S_1$ and $S_2$ are
invertible in $k[GL(m|n)]$. The fact that the determinant
$det(X_{11}-\Xi_{12}S_{22}(X_{22})\Xi_{21})$ is invertible is
contained in the proof of Proposition (3.2). The other determinant
can be proven invertible in the same way.
These determinants are also invertible in 
$k[SL(m|n)]$. In fact can be easily seen in $k[SL(m|n)]$
by observing that since $Ber=1$ and since:
$$
Ber=det(X_{11})det(X_{22}-\Xi_{21}S_{11}(X_{11})\Xi_{12})^{-1}
$$
(see [P]),
we have:
$$
\matrix{
det(X_{11}-\Xi_{12}S_{22}(X_{22})\Xi_{21})=det(S_{22}(X_{22}))^{-1}=
d^{m+1 \dots m+n}_{m+1 \dots m+n}
\cr\cr
det(X_{22}-\Xi_{21}S_{11}(X_{11})\Xi_{12})=
det((X_{11}))=d^{1 \dots m}_{1 \dots m}
}
$$
\medskip
{\bf Proposition (3.7)}.
{\it $k[GL(m|n)]$ and $k[SL(m|n)]$ are Hopf algebras with antipode $S$:}
$$
\matrix{
S\pmatrix{X_{11} & \Xi_{12} \cr \Xi_{21} & X_{22}}= 
\cr \cr\cr
\pmatrix{ S_1(X_{11}-\Xi_{12}S_{22}(X_{22})\Xi_{21}) &
         -S_{11}(X_{11})\Xi_{12}S_2(X_{22}-\Xi_{21}S_{11}(X_{11})\Xi_{12})
         \cr
         -S_{22}(X_{22})\Xi_{21}S_1(X_{11}-\Xi_{12}S_{22}(X_{22})\Xi_{21})
         & S_2(X_{22}-\Xi_{21}S_{11}(X_{11})\Xi_{12})}
}
$$
$$
S({d^{1...m}_{1...m}}^{-1})=d^{m+1...m+n}_{m+1...m+n} \qquad
S({d^{m+1...m+n}_{m+1...m+n}}^{-1})={d^{1...m}_{1...m}}
$$
{\bf Proof.}
One can check directly that this map is well defined and that
is an antipode (See [Be]).        
\medskip
{\bf Proposition (3.8)}.
{\it Let $A$ be a commutative superalgebra.

1. $Hom_{k-superalg}(k[GL(m|n)],A)$ is the group of automorphisms of 
$A^{m|n}$.

2. $Hom_{k-superalg}(k[SL(m|n)],A)$ is the group of automorphisms of 
$A^{m|n}$ with \break berezinian 1.

$k[GL(m|n)]$ and $k[SL(m|n)]$
are the representing objects for the functors $GL(m|n)$ 
and $SL(m|n)$ respectively.}

{\bf Proof}. Immediate.

\bigskip
{\bf 4. The quantum $GL(m|n)$}
\medskip

A quantum group is an Hopf algebra which is in general
neither commutative nor cocommutative.
According to this philosophy we can define a quantum supergroup in the 
same way.
\medskip
{\bf Definition (4.1)}. 
Let $A$ be a commutative (super)algebra over $k$.
A {\it formal deformation} of $A$ is a non commutative (super)algebra 
$A_q$ over $k_q=k[q,q^{-1}]$, $q$ being an (even) indeterminate,
such that $A_q/(q-1) \cong A$.
If $A$ in addition is an Hopf superalgebra
we will refer to such deformation as a {\it quantum (super)group}.
\medskip

We will call {\it quantum $GL(m|n)$} a formal 
deformation of the supercommutative
Hopf algebra $k[GL(m|n)]$ associated to  the affine
algebraic group $GL(m|n)$ (see \S 3). 
We cannot define quantum groups using the functor of
points as we did for supergroups. 
This is a consequence of the non commutativity of the
Hopf algebra associated  to the functor of points (see remark (2.2)(1)).

We will also show that it is possible 
to give a deformation of the Hopf algebra $k[GL(m|n)]$
in such a way that  natural coactions on quantum superspaces
are preserved.
 
In order to define the deformed algebra $k_q[GL(m|n)]$ we need first
to define the Manin matrix superalgebra $k_q[M(m|n)]$ 
introduced by Manin in [Ma3]
and to compute the commutation rules between the generators of 
$k_q[M(m|n)]$ and certain quantum determinants. 
\medskip
{\bf Definition (4.2).}
Define the following superalgebra ([Ma3]):
$$
k_q[M(m|n)]=_{def}k_q<x_{ij},\xi_{kl}>/I_M
$$
where $k_q<x_{ij},\xi_{kl}>$ denotes the free algebra over $k_q$
generated by the even variables
$x_{ij}$ for $1 \leq i,j \leq m$ or
                         $m+1 \leq i,j \leq m+n$ 
and by the odd variables
$\xi_{kl}$ for $1 \leq k \leq m, m+1 \leq l \leq m+n$ or
                         $m+1 \leq k \leq m+n, 1 \leq l \leq m$
satisfying the relation: $\xi_{kl}^2=0$.
 
Let's denote as before:
$$
a_{ij}=\cases{x_{ij} & if $1 \leq i,j \leq m$ or 
                         $m+1 \leq i,j \leq m+n$ \cr \cr
              \xi_{ij} & if $1 \leq i \leq m, m+1 \leq j \leq m+n$ or
                         $m+1 \leq i \leq m+n, 1 \leq j \leq m$}
$$
The ideal $I_M$ is generated by the relations ([Ma3]):
$$
\matrix{
a_{ij}a_{il}=(-1)^{\pi(a_{ij})\pi(a_{il})}
q^{(-1)^{p(i)+1}}a_{il}a_{ij}, \quad j < l \cr \cr
a_{ij}a_{kj}=(-1)^{\pi(a_{ij})\pi(a_{kj})}
q^{(-1)^{p(j)+1}}a_{kj}a_{ij}, \quad i < k \cr \cr
a_{ij}a_{kl}=(-1)^{\pi(a_{ij})\pi(a_{kl})}a_{kl}a_{ij}, \quad 
i< k,j > l \quad or \quad i > k,j < l \cr \cr
a_{ij}a_{kl}-(-1)^{\pi(a_{ij})\pi(a_{kl})}a_{kl}a_{ij}= \cr \cr
(-1)^{p(j)p(k)+p(j)p(l)+p(k)p(l)}(q^{-1}-(-1)^{\pi(a_{il})\pi(a_{jk})}q)
a_{jk}a_{il} \quad i<k,j<l
}
$$
where $p(i)=0$ if $1 \leq i \leq m$, $p(i)=1$ otherwise
and $\pi(a_{ij})$ denotes the parity of $a_{ij}$.
 
Notice that for $q=1$ this gives us the superalgebra defined
in the example (2.3)(2) representing the functor $M(m|n)$. 
\medskip
{\bf Definition (4.3).}
Define the {\it quantum superspace} $k_q^{m|n}$ as the ring generated over
$k_q$ by the even elements $x_1 \dots x_m$ and the odd elements
$\xi_1 \dots \xi_n$ subject to the relations [Ma3]:
$$
\matrix{
x_ix_j-q^{-1}x_jx_i, \qquad 1 \leq i < j \leq  m, \cr \cr
x_i\xi_k-q^{-1}\xi_kx_i, \qquad 1 \leq i \leq m, \quad 
m+1 \leq k \leq n+m, \cr \cr
\xi_k^2, \qquad \xi_k\xi_l+q^{-1}\xi_l\xi_k, \qquad m+1 \leq k < l \leq  m+n.
}
$$
Define also {\it dual quantum superspace} $(k_q^{m|n})^*$  
as the ring generated over
$k_q$ by the even elements $y_1 \dots y_n$ and the  odd elements
$\eta_1 \dots \eta_m$ subject to the relations ([Ma3]):
$$
\matrix{
y_iy_j-qy_jy_i, \qquad 1 \leq i < j \leq  m, \cr \cr
y_i\eta_k-q\eta_ky_i, \qquad 1 \leq i \leq m, \quad
m+1 \leq k \leq n+m, \cr \cr
\eta_k^2, \qquad \eta_k\eta_l+q\eta_l\eta_k, \qquad m+1 \leq k < l \leq  m+n.
}
$$
\vskip1ex
{\bf Observation (4.4).}
The superalgebra $k_q[M(m|n)]$ admits a bialgebra structure with
comultiplication and counit given by:
$$
\matrix{
\D(a_{ij})=\sum a_{ik} \otimes a_{kj} &
\e(a_{ij})=\d_{ij}
}
$$
and a coaction on the quantum spaces $k_q^{m|n}$ and $(k_q^{m|n})^*$
([Ma3]).
\vskip1ex
{\bf Observation (4.5).}
Let's examine some of the relations that generate the ideal $I_M$.

For $1 \leq i,j,k,l \leq m$:
$$
\matrix{
x_{ij}x_{il}=q^{-1}x_{il}x_{ij}, \quad j<l, \qquad
x_{ij}x_{kj}=q^{-1}x_{kj}x_{ij}, \quad k<i, \cr\cr
x_{ij}x_{kl}=x_{kl}x_{ij}, \quad i<k,j>l \quad \hbox{or} \quad i>k,j<l, \cr \cr
x_{ij}x_{kl}-x_{kl}x_{ij}=(q^{-1}-q)x_{il}x_{kj}, \quad i<k,j<l. 
}
$$
which are the usual Manin relations (see [Ma2]). We denote
the two-sided ideal generated by them as $I_M^{1 \dots m}(q)$.

For $m+1 \leq i,j,k,l \leq m+n$:
$$
\matrix{
x_{ij}x_{il}=qx_{il}x_{ij}, \quad j<l \qquad
x_{ij}x_{kj}=qx_{kj}x_{ij}, \quad k<i \cr\cr
x_{ij}x_{kl}=x_{kl}x_{ij}, \quad i<k,j>l \quad \hbox{or} \quad i>k,j<l \cr \cr
x_{ij}x_{kl}-x_{kl}x_{ij}=(q^{-1}-q)x_{il}x_{kj}, \quad i<k,j<l.
}
$$
These relations are the usual Manin relations 
where $q$ is replaced with $q^{-1}$.
We denote the ideal generated by them as $I_M^{m+1 \dots m+n}(q^{-1})$.

Let's define:
$$
\matrix{
D_{1 \dots m}^{1 \dots m}=_{def}\sum_{\s \in S_m}(-q)^{-l(\s)}
x_{1\s(1)} \dots x_{m\s(m)} \cr \cr
D_{m+1 \dots m+n}^{m+1 \dots m+n}=_{def}\sum_{\s \in S_n}(-q)^{l(\s)}
x_{m+1,m+\s(1)} \dots x_{m+n,m+\s(n)}. 
}
$$

$D_{1 \dots m}^{1 \dots m}$ and $D_{m+1 \dots m+n}^{m+1 \dots m+n}$
represent respectively the quantum determinants of the quantum matrix
bialgebras:
$$
\matrix{
M_q(1,m)=k_q<x_{ij}>/I_M^{1 \dots m}(q), \qquad 1 \leq i,j \leq m,
\cr\cr\cr
M_{q^{-1}}(m+1,m+n)=k_q<x_{kl}>/I_M^{m+1 \dots m+n}(q^{-1}), \qquad 
m+1 \leq k,l \leq m+n.
}
$$ 
$D_{1 \dots m}^{1 \dots m}$ and $D_{m+1 \dots m+n}^{m+1 \dots m+n}$
are central elements in $M_q(1,m)$ and $M_{q^{-1}}(m+1,m+n)$
respectively ([PW] pg 50).

Further properties of these algebras and their quantum
determinants are studied in [PW] ch. 5.

We define:
$$
\matrix{
GL_{q}(1,m)=M_q(1,m)<{D_{1 \dots m}^{1 \dots m}}^{-1}>/
(D_{1 \dots m}^{1 \dots m}
{D_{1 \dots m}^{m \dots m}}^{-1}-1), 
}
$$
$$
GL_{q^{-1}}(m+1,m+n)={M_{q^{-1}}(m+1,m+n)
<{D_{m+1 \dots m+n}^{m+1 \dots m+n}}^{-1}>
\over
(D_{m+1 \dots m+n}^{m+1 \dots m+n}
{D_{m+1 \dots m+n}^{m+1 \dots m+n}}^{-1}-1).
}
$$
These are Hopf algebras, their antipodes $S_{11}^q$,
$S_{22}^q$ are explicitly calculated
in [PW] pg 57.

\vskip1ex
{\bf Definition (4.6)}.
{\it Quantum general linear supergroup}.
$$
k_q[GL(m|n)]=_{def}
k_q[M(m|n)]
<{D_{1 \dots m}^{1 \dots m}}^{-1},{D_{m+1 \dots m+n}^{m+1 \dots m+n}}^{-1}>
$$
where ${D_{1 \dots m}^{1 \dots m}}^{-1}$, 
${D_{m+1 \dots m+n}^{m+1 \dots m+n}}^{-1}$ are even indeterminates such
that:
$$
\matrix{
{D_{1 \dots m}^{1 \dots m}}{D_{1 \dots m}^{1 \dots m}}^{-1}=1=
{D_{1 \dots m}^{1 \dots m}}^{-1}{D_{1 \dots m}^{1 \dots m}},\cr\cr\cr
{D_{m+1 \dots m+n}^{m+1 \dots m+n}}
{D_{m+1 \dots m+n}^{m+1 \dots m+n}}^{-1}=1=
{D_{m+1 \dots m+n}^{m+1 \dots m+n}}^{-1}
{D_{m+1 \dots m+n}^{m+1 \dots m+n}}
}
$$
where $X_{11}=(x_{ij})_{1 \leq i,j \leq m}$,
$X_{22}=(x_{ij})_{m+1 \leq i,j \leq m+n}$ are matrices of even
indeterminates and
$\Xi_{12}=(\xi_{kl})_{1 \leq k \leq m, m+1 \leq l \leq m+n}$,
$\Xi_{12}=(\xi_{kl})_{m+1 \leq k \leq m+n, 1 \leq l \leq m}$
are matrices of odd indeterminates.
The $x_{ij}$'s, $\xi_{kl}$'s are the generators of
$k_q[M(m|n)]$.

$det_q(M)$ denotes the determinant of the Manin matrix $M$ (that is
a matrix of indeterminates satisfying the Manin commutation relations). 
One can verify ([PW] ch. 5) that $det_q(S_{22}^q(X_{22}))=
{D_{m+1 \dots m+n}^{m+1 \dots m+n}}^{-1}$. 

By an abuse of language we will denote with the same
letters $x_{ij}$, $\xi_{kl}$ and $a_{ij}$ the indeterminates
generators of the rings $k_q[M(m|n)]$ and $k_q[GL(m|n)]$ the context
making clear where the elements sit.
\medskip
We also define:
$$
Ber_q=_{def}
det_q(S_{22}^q(X_{22}))det_{q}(X_{11}-\Xi_{12}S_{22}^q(X_{22})\Xi_{21}).
$$
called the {\it quantum berezinian}. $Ber_q$ is invertible in
$k_q[GL(m|n)]$, the proof is a small variation of the one in (3.2).

\medskip
We now want to explicitly give the coalgebra structure for the
ring $k_q[GL(m|n)]$. We need some lemmas on commutation relations.
\medskip
%
%

{\bf Lemma (4.7)}. {\it
$$
\matrix{
1) \quad
D^{1...m}_{1...m}\xi_{ij}=q^{-1}\xi_{ij}D^{1...m}_{1...m}, \qquad
{D^{1...m}_{1...m}}^{-1}\xi_{ij}=q\xi_{ij}{D^{1...m}_{1...m}}^{-1},
\cr\cr
2) \quad 
D^{m+1...m+n}_{m+1...m+n}\xi_{ij}=q^{-1}
\xi_{ij}D^{m+1...m+n}_{m+1...m+n}, \qquad
{D^{m+1...m+n}_{m+1...m+n}}^{-1}\xi_{ij}=
q\xi_{ij}{D^{m+1...m+n}_{m+1...m+n}}^{-1},
}
$$
where 
$1 \leq i \leq m$, $m+1 \leq j \leq m+n$ or
$m+1 \leq i \leq m+n$, $1 \leq j \leq m$.
}

{\bf Proof}. Let's prove (1).
It is enough to prove the first commutation relation.
We will first do for $1 \leq i \leq m$, $m+1 \leq j \leq m+n$.

By induction on $m$. For $m=1$ it is a direct simple check.
Assume for now that $1 \leq i < m$, $m+1 \leq j \leq m+n$ (the case
$i=m$ is treated separately).

Let $D^{1 \dots \hat \alpha \dots m}_{1 \dots \hat \beta \dots m}$
denote the quantum determinant of the quantum minor obtained
from the quantum matrix $X_{11}$ by suppressing row $\alpha$
and column $\beta$.

By [PW] formula at pg 47 on Laplace expansion of quantum 
determinants we have:
$$
\matrix{
D^{1...m}_{1...m}\xi_{ij}=
\sum_{s=1}^m(-q)^{s-m}D^{1 \dots \hat s \dots m}_{1 \dots m-1}x_{ms}\xi_{ij}=
\cr\cr
=\sum_{s=1}^m(-q)^{s-m}D^{1 \dots \hat s \dots m}_{1 \dots m-1}\xi_{ij}x_{ms}.
}
$$
By induction 
$$
D^{1 \dots \hat s \dots m}_{1 \dots m-1}\xi_{ij}=
q^{-1}\xi_{ij}D^{1 \dots \hat s \dots m}_{1 \dots m-1}
$$
hence we have our result. 

Now the case 
$i=m$, $m+1 \leq j \leq m+n$. By [PW] formula at pg 47 we have:
$$
\matrix{
D^{1...m}_{1...m}\xi_{mj}=
\sum_{s=1}^m(-q)^{s-1}D^{1 \dots \hat s \dots m}_{2 \dots m}
x_{1s}\xi_{mj}= \cr\cr
=\sum_{s=1}^m(-q)^{s-1}D^{1 \dots \hat s \dots m}_{2 \dots m}
[\xi_{mj}x_{1s}+(q^{-1}-q)\xi_{1j}x_{ms}].
}
$$
By induction: 
$$
D^{1 \dots \hat s \dots m}_{2 \dots m}\xi_{mj}=
q^{-1}\xi_{mj}D^{1 \dots \hat s \dots m}_{2 \dots m}.
$$
So we have:
$$
\matrix{
D^{1...m}_{1...m}\xi_{mj}=q^{-1}\xi_{mj}D^{1...m}_{1...m}+
(q^{-1}-q)\sum_{s=1}^m(-q)^{s-1}D^{1 \dots \hat s \dots m}_{2 \dots m}
x_{ms}\xi_{1j}].
}
$$
But notice that ([PW] pg 47):
$$
\sum_{s=1}^m(-q)^{s-1}D^{1 \dots \hat s \dots m}_{2 \dots m}
x_{ms}=\d_{1m}D^{1...m}_{1...m}=0
$$
hence we have our result.

The case $m+1 \leq i \leq m+n$, $1 \leq j \leq m$ is done in a
similar way.

The proof of (2) goes along the same lines.
\medskip
{\bf Lemma (4.8)}.
{\it
$$
\matrix{
1.a)  \quad D^{1 \dots m}_{1 \dots m}
a_{1,k_1} \dots a_{m,k_m}=
q^{-t}a_{1,k_1} \dots a_{m,k_m}D^{1 \dots m}_{1 \dots m} 
\cr\cr\cr
1.b) \quad D^{1 \dots m}_{1 \dots m}a_{k_1,1} \dots a_{k_m,m}=
q^{-t}a_{k_1,1} \dots a_{k_m,m}D^{1 \dots m}_{1 \dots m}
\cr\cr\cr
2.a) \quad D^{m+1 \dots m+n}_{m+1 \dots m+n}
a_{m+1,l_1} \dots a_{m+n,l_n}=
q^{-s}a_{m+1,l_1} \dots a_{m+n,l_n}
D^{m+1 \dots m+n}_{m+1 \dots m+n}
\cr\cr\cr
2.b) \quad D^{m+1 \dots m+n}_{m+1 \dots m+n}a_{l_1,m+1} \dots 
a_{l_n,m+n} =
q^{-s}a_{l_1,m+1} \dots a_{l_n,m+n} 
D^{m+1 \dots m+n}_{m+1 \dots m+n}
\cr\cr
}
$$
where $t$ is such that $m-t=|\{k_1 \dots k_m\} \cap \{1 \dots m\}|$, 
$1 \leq k_1 < \dots < k_m \leq m+n$ and
$s$ is such that $n-s=|\{l_1 \dots l_n\} \cap \{m+1 \dots m+n\}|$,
$1 \leq l_1 < \dots < l_n \leq m+n$.
}

{\bf Proof}. Let's prove (1.a). 
By lemma (4.7) since $D^{1 \dots m}_{1 \dots m}$ commutes
with $x_{ij}$ ([PW] pg 50), 
$1 \leq i,j \leq m$ we have that for $1 \leq i \leq m$:
$$
D^{1 \dots m}_{1 \dots m}a_{i,k_j}=
\cases{a_{i,k_j} D^{1 \dots m}_{1 \dots m} & for $1 \leq k_j \leq m$ \cr \cr
       q^{-1}a_{i,k_j} D^{1 \dots m}_{1 \dots m} & otherwise}
$$
Hence we have the result. The proofs of (1.b) and (2.a), (2.b) are the same.
\medskip
As in the supercommutative case we can write:
$$
\D(D^{1 \dots m}_{1 \dots m})=
\sum_{i=1}^m R_i
$$
with 
$$
\matrix{
R_i=\sum_{\s \in S_m,(k_1, \dots , k_m)\in \rho_i}
(-q)^{-l(\s)}a_{1,k_1} \dots a_{m,k_m} \otimes
a_{k_1,\s(1)} \dots a_{k_m,\s(m)} \cr\cr
\rho_i=\{(k_1, \dots , k_m) | 
1 \leq k_1, \dots , k_m \leq m+n, \quad
|\{k_1, \dots , k_m\} \cap \{1 \dots m\}|=m-i\}
}
$$
Similarly one can write:
$$
\D(D_{m+1 \dots m+n}^{m+1 \dots m+n})=\sum_{i=1}^n S_i
$$
with
$$
\matrix{
s_i=\sum_{\t \in S_n,(l_1, \dots , l_n)\in \sigma_i}
(-q)^{l(\t)}a_{1,l_1} \dots a_{n,l_n} \otimes
a_{l_1,\t(1)} \dots a_{l_n,\t(n)} \cr\cr
\sigma_i=\{(l_1, \dots , l_n)| 1 \leq l_1, \dots , l_n \leq m+n,
\quad
|\{l_1, \dots , l_n\} \cap \{m+1 \dots m+n\}|=n-i\}.
}
$$

\medskip
{\bf Lemma (4.9)}.{\it
 
1) $R_0R_i=q^{-2i}R_iR_0, \qquad R_0^{-1}R_i=q^{2i}R_iR_0^{-1}$

2) $S_0S_i=q^{-2i}S_iS_0, \qquad S_0^{-1}S_i=q^{2i}S_iS_0^{-1}$
}

{\bf Proof}. Immediate from lemma (4.8) noting that
$R_0=D^{1 \dots m}_{1 \dots m}\otimes D^{1 \dots m}_{1 \dots m}$
and $S_0=D^{m+1 \dots m+n}_{m+1 \dots m+n}\otimes 
D^{m+1 \dots m+n}_{m+1 \dots m+n}$.
\medskip
{\bf Lemma (4.10)}. 

{\it 
1. $R_0(R_1+ \dots +R_m)=(q^{-2}R_1+\dots+q^{-2m}R_m)R_0$

2. $(R_1+ \dots +R_m)R_0^{-1}=R_0^{-1}(q^{-2}R_1+\dots+q^{-2m}R_m)$

3. $S_0(S_1+ \dots +S_n)=(q^{-2}S_1+\dots+q^{-2n}S_n)S_0$

4. $(S_1+ \dots +S_n)S_0^{-1}=S_0^{-1}(q^{-2}S_1+\dots+q^{-2n}S_n)$
}

{\bf Proof}.
This is an immediate application of lemma (4.9).
\medskip
{\bf Lemma (4.11)}. {\it Let $a_{i_1^s,j_1^s} \dots a_{i_m^s,j_m^s} \in
k_q[GL(m|n)]$, for all $s$, $1 \leq s \leq 2mn+1$. Assume
that for all $s$
at least one $(i_k^s,j_k^s)$ is such that
$1 \leq i_k \leq m$, $m \leq j_k \leq m+n$ or $m+1 \leq i_k \leq m+n$,
$1 \leq j_k \leq m$, that is $a_{i_k,j_k}^s$ is odd. Then
$$
\prod_{s=1}^{2mn+1}a_{i_1^s,j_1^s} \dots a_{i_m^s,j_m^s} =0
$$
}

{\bf Proof}. Let  $i$ be the map:
$$
\matrix{
k_q[M(m|n)] & \mapemb{i} &
k_q[M(m|n)]<{D^{1...m}_{1...m}}^{-1},{D^{m+1...m+n}_{m+1...m+n}}^{-1}>
}
$$

Let $X \in k_q[GL(m|n)]$ such that there exists a $X_0 \in k_q[M(m|n)]$, 
$i(X_0)=X$. 
If one wants
to prove $X=0$ it is enough to show $X_0=0$. 

Now let 
$$
X= \prod_{s=1}^{2mn+1}a_{i_1^s,j_1^s} \dots a_{i_m^s,j_m^s}.
$$
If $X_0= \prod_{s=1}^{2mn+1}a_{i_1^s,j_1^s} \dots a_{i_m^s,j_m^s} 
\in k_q[M(m|n)]$ we have that $i(X_0)=X$. 

By the previous argument it is enough to show that 
$X_0=0$.

By [Ma3] pg 172 we have that 
$\{1,a_{i_1,j_1} \dots a_{i_r,j_r} \}_{\{(i_1,j_1) 
\leq \dots \leq (i_r,j_r), r \geq 1\}}$
form a basis for \break $k_q[M(m|n)]$, where $\leq$ is 
a suitable ordering on the indeces $(i_k,j_k)$.

Since the Manin relations are homogeneous we have:
$$
X_0=\sum_{(r_1,s_1) \leq \dots \leq (r_t,s_t),
c_{(r_1,s_1) \dots (r_t,s_t)} \in k_q}
c_{(r_1,s_1) \dots (r_t,s_t)} a_{r_1,s_1} \dots a_{r_t,s_t} \qquad (\star)
$$
where the sum is taken over a suitable set of 
$((r_1,s_1) \dots (r_t,s_t))$ with $t=m(2mn+1)$.

By hypothesis each term in the sum $(\star)$ contains 
$2mn+1$ odd indeterminates and since
the indeces are ordered, it is 0.
\medskip 
{\bf Proposition (4.12)}. {\it Coalgebra structure for $k_q[GL(m|n)]$.}

{\it $k_q[GL(m|n)]$ is a coalgebra with comultiplication:
$$
\matrix{
\D(a_{ij})=\sum a_{ik} \otimes a_{kl} \cr\cr\cr
\D({D^{1...m}_{1...m}}^{-1})=
\sum_{i=1}^{2mn+2} (-1)^{i-1} R_0^{-1}
[R_0^{-1}(q^{-2}R_1+\dots +q^{-2m}R_m)]^{i-1}
\cr \cr \cr
\D({D^{m+1...m+n}_{m+1...m+n}}^{-1})=
\sum_{i=1}^{2mn+2} (-1)^{i-1} S_0^{-1}[S_0^{-1}
(q^{-2}S_1+\dots +q^{-2n}S_n)]^{i-1}
}
$$
and counit:
$$
\matrix{
\e(a_{ij})=\d_{ij}, \qquad
\e({D^{1...m}_{1...m}}^{-1})=1, \qquad \e({D^{m+1...m+n}_{m+1...m+n}}^{-1})=1
}
$$
}

{\bf Proof.}
$\D$ is well defined on all the commutation relations among the generators.

We only need to check that 
$$
\matrix{
\D({D_{1 \dots m}^{1 \dots m}})\D({D_{1 \dots m}^{1 \dots m}}^{-1})=1=
\D({D_{1 \dots m}^{1 \dots m}}^{-1})\D({D_{1 \dots m}^{1 \dots m}}),
\cr\cr\cr
\D({D_{m+1 \dots m+n}^{m+1 \dots m+n}})
\D({D_{m+1 \dots m+n}^{m+1 \dots m+n}}^{-1})=1=
\D({D_{m+1 \dots m+n}^{m+1 \dots m+n}}^{-1})
\D({D_{m+1 \dots m+n}^{m+1 \dots m+n}})
}
$$
Let's check the first one.
$$
\matrix{
\D({D_{1 \dots m}^{1 \dots m}})
\D({D_{1 \dots m}^{1 \dots m}}^{-1})=
}
$$
$$
[R_0+(R_1+\dots+ R_m)][R_0^{-1}-R_0^{-2}(q^{-2}R_1 +\dots +q^{-2m}R_m)+ \dots]=
$$
$$
=1 \otimes 1+(R_1+\dots R_m)R_0^{-1}-
R_0^{-1}(q^{-2}R_1 +\dots +q^{-2m}R_m)+
$$
$$
-(R_1+\dots+ R_m)R_0^{-2}(q^{-2}R_1 +\dots +q^{-2m}R_m)+ \dots
$$
We obtain a telescopic sum. In fact, let's see the generic terms:
$$
\matrix{
G=[R_0+(R_1+\dots+ R_m)][\dots +(-1)^{i-1} R_0^{-1}
[R_0^{-1}(q^{-2}R_1+\dots +q^{-2m}R_m)]^{i-1}+
\cr\cr\cr
+(-1)^{i} R_0^{-1}[R_0^{-1}(q^{-2}R_1+\dots +q^{-2m}R_m)]^{i}+\dots]=
\cr\cr\cr
=\dots +(-1)^{i-1}[R_0^{-1}(q^{-2}R_1+\dots +q^{-2m}R_m)]^{i-1}+
\cr\cr\cr
+(-1)^{i-1}(R_1+\dots+ R_m)R_0^{-1}
[R_0^{-1}(q^{-2}R_1+\dots +q^{-2m}R_m)]^{i-1}+
\cr\cr\cr
+(-1)^{i}[R_0^{-1}(q^{-2}R_1+\dots +q^{-2m}R_m)]^{i}+
\cr\cr\cr
+(-1)^{i}(R_1+\dots+ R_m)R_0^{-1}[R_0^{-1}(q^{-2}R_1+\dots +q^{-2m}R_m)]^{i}.
}
$$
By lemma (4.10)(2) we have that
$(R_1+\dots+ R_m)R_0^{-1}=
R_0^{-1}(q^{-2}R_1+\dots +q^{-2m}R_m)$
so the second and third term in the sum $G$ cancel each other.

The last term of the sum is 0 by lemma (4.11).

This completes the check for the first part of the first relation.

Let's see the second part.
$$
\matrix{
\D({D_{1 \dots m}^{1 \dots m}}^{-1})
\D({D_{1 \dots m}^{1 \dots m}})=\cr\cr\cr
=
[R_0^{-1}-R_0^{-2}(q^{-2}R_1 +\dots +q^{-2m}R_m)+ \dots][R_0+(R_1+\dots+ R_m)]=
\cr\cr\cr
=1 \otimes 1+R_0^{-1}(R_1+\dots R_m)-
R_0^{-2}(q^{-2}R_1 +\dots +q^{-2m}R_m)R_0+
\cr\cr\cr
-R_0^{-2}(q^{-2}R_1 +\dots +q^{-2m}R_m)(R_1+\dots+ R_m)+ \dots
}
$$
By lemma (4.10)(1) the second and third term cancel out.
Let's see as before the generic terms $G'$:
$$
\matrix{
G'=[\dots +(-1)^{i-1} R_0^{-1}
[R_0^{-1}(q^{-2}R_1+\dots +q^{-2m}R_m)]^{i-1}+
\cr\cr\cr
+(-1)^{i} R_0^{-1}[R_0^{-1}(q^{-2}R_1+\dots +q^{-2m}R_m)]^{i}+\dots]
[R_0+(R_1+\dots+ R_m)]=
\cr\cr\cr
=\dots+(-1)^{i-1} R_0^{-1}
[R_0^{-1}(q^{-2}R_1+\dots +q^{-2m}R_m)]^{i-1}R_0+
\cr\cr\cr
+(-1)^{i-1} R_0^{-1}
[R_0^{-1}(q^{-2}R_1+\dots +q^{-2m}R_m)]^{i-1}
(R_1+\dots+ R_m)+
\cr\cr\cr
+(-1)^{i} R_0^{-1}[R_0^{-1}(q^{-2}R_1+\dots +q^{-2m}R_m)]^{i}R_0+
\cr\cr\cr
+(-1)^{i} R_0^{-1}[R_0^{-1}(q^{-2}R_1+\dots +q^{-2m}R_m)]^{i}
(R_1+\dots+ R_m)+\dots
}
$$
We now look at the second and third term in $G'$.
$$
\matrix{
+(-1)^{i-1}R_0^{-1}
[R_0^{-1}(q^{-2}R_1+\dots +q^{-2m}R_m)]^{i-1}
(R_1+\dots+ R_m)+
\cr\cr\cr
+(-1)^{i}R_0^{-1}[R_0^{-1}(q^{-2}R_1+\dots +q^{-2m}R_m)]^{i-1}
R_0^{-1}(q^{-2}R_1+\dots +q^{-2m}R_m)R_0.
}
$$
By lemma (4.10)(1) we have that $(q^{-2}R_1+\dots +q^{-2m}R_m)R_0=
R_0(R_1+\dots+ R_m)$, so the second and third term cancel each
other. The last term of the sum is 0 by lemma (4.11).

The second relation can be checked in the same way.
\medskip
{\bf Proposition (4.13)}. {\it $k_q[GL(m|n)]$ admits a coaction on
$k^{m|n}_q$, ${k^{m|n}_q}^*$.}

{\bf Proof}. Immediate.
\medskip
Now we need to give the antipode on $k_q[GL(m|n)]$.
\medskip
{\bf Proposition (4.14)}. {\it $k_q[GL(m|n)]$ is an Hopf algebra 
with antipode $S^q$ given by:
$$
\matrix{
S^q(X_{11})=
S_{11}^q(X_{11})+S_{11}^q(X_{11})\Xi_{12}S_2^q(X_{22}-\Xi_{21}
S_{11}^q(X_{11})\Xi_{12})\Xi_{21}S_{11}^q(X_{11}) \cr\cr\cr
S^q(\Xi_{12})=-S_{11}^q(X_{11})\Xi_{12}S_2^q(X_{22}-\Xi_{21}
S_{11}^q(X_{11})\Xi_{12}) \cr\cr\cr
}
$$
$$
\matrix{
S^q(\Xi_{21})=-S_2^q(X_{22}-\Xi_{21}
S_{11}^q(X_{11})\Xi_{12})\Xi_{21}S_{11}^q(X_{11})\cr\cr\cr
S^q(X_{22})=S_2^q(X_{22}-\Xi_{21}
S_{11}^q(X_{11})\Xi_{12})\cr\cr\cr
}
$$
$$
\matrix{
S^q({D_{1 \dots m}^{1 \dots m}}^{-1})=
{D_{m+1 \dots m+n}^{m+1 \dots m+n}}\cr\cr\cr
S^q({D_{m+1 \dots m+n}^{m+1 \dots m+n}}^{-1})=
{D_{1 \dots m}^{1 \dots m}}
}
$$
where $X_{22}-\Xi_{21}S_{11}^q(X_{11})\Xi_{12}$ is a quantum
matrix (see [P] \S 4) and $S^q_2$ denotes its quantum antipode.}

{\bf Proof}. See [P] \S 4.
\medskip
{\bf Remark (4.15)}. In the ring $k_q[GL(m|n)]$ the parameter
$q$ can be specialized to any value in $k^{\times}$.
Hence in the definition of $k_q[GL(m|n)]$ $q$ can be also
taken as any element in $k^{\times}$.
\medskip
{\bf Example (4.16)}: $k_q[GL(1|1)]$.
$$
k_q[GL(1|1)]=k<x_{11},\xi_{12},\xi_{21}, x_{22},x_{11}^{-1},
x_{22}^{-1}>/I
$$
where $I$ is the ideal generated by the relations:
$$
\matrix{
x_{11}\xi_{12}=q^{-1}\xi_{12}x_{11}, \qquad
x_{11}\xi_{21}=q^{-1}\xi_{21}x_{11}, \cr\cr
\xi_{21}x_{22}=qx_{22}\xi_{21},\qquad
\xi_{12}x_{22}=qx_{22}\xi_{12},\cr\cr
x_{11}x_{22}-x_{22}x_{11}=(q-q^{-1})\xi_{12}\xi_{21}\cr\cr
\xi_{12}\xi_{21}=-\xi_{21}\xi_{12},\qquad
\xi_{12}^2=\xi_{21}^2=0.\cr\cr
}
$$
The coalgebra structure is the following.
$$
\matrix{
\D(a_{ij})=\sum_k a_{ik} \otimes a_{kj}\cr\cr
\D(x_{11}^{-1})=
x_{11}^{-1}\otimes x_{11}^{-1}-q^{-2}x_{11}^{-2}\xi_{12}\otimes
x_{11}^{-2}\xi_{21}, \cr\cr
\D(x_{22}^{-1})=
x_{22}^{-1}\otimes x_{22}^{-1}-q^2x_{22}^{-2}\xi_{21}\otimes
x_{22}^{-2}\xi_{12}, \cr\cr
\e(x_{11}^{-1})=\e(x_{22}^{-1})=1,
\qquad
\e\pmatrix{x_{11} & \xi_{12} \cr
\xi_{21} & x_{22}}=
\pmatrix{1 & 0 \cr
0 & 1},
}
$$
$$
S^q\pmatrix{
x_{11} & \xi_{12} \cr\cr
\xi_{21} & x_{22}}=
\pmatrix{ (x_{11}-\xi_{12}x_{22}^{-1}\xi_{21})^{-1} &
         -x_{11}^{-1}\xi_{12}(x_{22}-\xi_{21}x_{11}^{-1}\xi_{12})^{-1}
         \cr\cr
         -x_{22}^{-1}\xi_{21}(x_{11}-\xi_{12}x_{22}^{-1}\xi_{21})^{-1}
         & (x_{22}-\xi_{21}x_{11}^{-1}\xi_{12})^{-1}}
$$
$$
\matrix{
S(x_{11}^{-1})=x_{22} \cr\cr
S(x_{22}^{-1})=x_{11}
}
$$

$k_q[GL(1|1)]$ admits a coaction on the quantum spaces:
$$
\matrix{
k_q^{1|1}=k<x,\xi>/(x\xi-q^{-1}\xi x) \cr\cr
(k_q^{1|1})^*=k<y,\eta>/(y\eta-q^{-1}\eta y)
}
$$
\medskip
In this particular case is it immediate to construct also a deformation
for the Hopf algebra $k[SL(1|1)]$. In fact since the
quantum berezinian:
$$
Ber_q=x_{22}^{-1}(x_{11}-\xi_{12}x_{22}^{-1}\xi_{21})
$$
is a central element in $k_q[M(m|n)]$ and in $k_q[GL(m|n)]$
we can define:
$$
k_q[SL(1|1)]=
k<x_{11},\xi_{12},\xi_{21}, x_{22},x_{11}^{-1},
x_{22}^{-1}>/I'
$$
where $I'$ is the two-sided ideal generated by the same relations
as $I$ together with the extra relation $Ber_q=1$.

The comultiplication, counit, antipode are the same as $k_q[GL(1|1)]$,
in fact one can check directly that they are still well defined.

\medskip
We plan to construct in a forthcoming paper the deformation
$k_q[SL(m|n)]$ and its relation with $k_q[GL(m|n)]$.
\bigskip
\centerline{\bf REFERENCES}
\bigskip

\item{[Be]} Berezin, F. A., {\it Introduction to superanalysis},
            Kluwer Academic Publisher (1987).

\item{[Bo]} Borel, A., {\it Linear algebraic groups}, GTM Springer-Verlag,
            (1991).
 
\item{[CP]} Chari, V., Pressley, A., {\it A guide to quantum groups},
            Cambridge University Press, (1994).

\item{[DM]} Deligne, P., Morgan, J.
		{\it Supersymmetry}, Quantum Fields and strings. A 
            course for mathematicians, Vol 1, AMS, (1999).

\item{[Fi]} Fioresi, R., {\it Supergroups, quantum supergroups and their
        homogeneous spaces}, Modern Physics Letters A, 16, no. 4-6,
        269-274, (2001).

\item{[Ha]} Hartshorne, R., {\it Algebraic geometry}, GTM Springer-Verlag,
            (1977).

\item{[Ko]} Kostant, B., {\it Graded manifolds, graded Lie theory
	    and prequantization}, Springer LNM 570, (1977).

\item{[Le]} Leites, D. A., {\it Introduction to the theory of supermanifolds},
            Russian Math. Surveys, {\bf 35}, no. 1, 1-64, (1980).

\item{[LS]} Lyubashenko V., Sudbery A., {\it Quantum supergroups of 
            $GL(m|n)$ type}, Duke Math. J., {\bf 90}, no. 1, (1997).  

\item{[Ma1]} Manin, Yu., {\it Gauge field theory and complex geometry},
             Springer Verlag, (1988).

\item{[Ma2]} Manin, Yu., {\it Topics in non commutative geometry},
             Princeton University Press, (1991). 

\item{[Ma3]} Manin, Yu., {\it Multiparametric quantum deformation of
             the general linear supergroup}, Comm. Math. Phy., {\bf 123}, 
             163-175, (1989).

\item{[Ma4]} Manin, Yu., {\it Quantum groups and non commutative geometry}, 
             Centre de Reserches Mathematiques Montreal, 49, (1988). 
        
\item{[Mo]} Montgomery, S. {\it Hopf algebras and their actions on rings},
	AMS, Providence, (1993).
      
\item{[P]} Phung, Ho Hai, {\it On the structure of the quantum supergroups
            $GL_q(m|n)$}, q-alg$\slash$9511-23, (1999).

\item{[PW]} Parshall, B.; Wang, J. P., {\it Quantum linear groups}, 
Mem. Amer. Math. Soc. 89, no. 439, (1991).

\end